\documentclass[letterpaper, 10 pt, conference]{ieeeconf}  

\IEEEoverridecommandlockouts                              
\usepackage{graphics} 
\usepackage{epsfig} 
\usepackage{mathptmx} 
\usepackage{times} 
\usepackage{amsmath} 
\usepackage{amssymb}  
\usepackage{algorithmic}
\usepackage[nolist]{acronym}
\usepackage{cite}
\usepackage[bookmarks=true]{hyperref}
\usepackage{xcolor}
\usepackage{tabularx}
\usepackage[ruled,linesnumbered]{algorithm2e}
\usepackage{url}
\usepackage{times}
\usepackage{amsfonts}
\usepackage{mathrsfs}
\usepackage{subfigure}
\usepackage{wrapfig}
\usepackage{comment}
\usepackage{xspace}
\usepackage{enumerate}
\usepackage{epsfig}
\usepackage{epstopdf}
\usepackage{booktabs}
\usepackage{multicol}
\usepackage{graphicx}
\usepackage{caption}
\usepackage{comment}

\newtheorem{theorem}{\bf{Theorem}}

\DeclareMathOperator{\tr}{\rm{tr}}
\newcommand{\rd}{{\mathrm d}}
\newcommand{\rT}{{\mathrm T}}
\newcommand{\E}{\mathbb{E}}
\newcommand{\R}{\mathbb{R}}
\newcommand{\calU}{{\cal{U}}}
\newcommand{\calV}{{\cal{V}}}

\begin{acronym}
\acro{DNN}{Deep Neural Network}
\acro{DP}{Dynamic Programming}
\acro{FBSDE}{Forward-Backward Stochastic Differential Equation}
\acro{LSTM}{Long-Short Term Memory}
\acro{DDP}{Differential Dynamic Programming}
\acro{HJB}{Hamilton-Jacobi-Bellman}
\acro{HJI}{Hamilton-Jacobi-Isaacs}
\acro{PDE}{Partial Differential Equation}
\acro{PI}{Path Integral}
\acro{NN}{Neural Network}
\acro{GPs}{Gaussian Processes}
\acro{SOC}{Stochastic Optimal Control}
\acro{RL}{Reinforcement Learning}
\acro{MPOC}{Model Predictive Optimal Control}
\acro{IL}{Imitation Learning}
\acro{RNN}{Recurrent Neural Network}
\acro{DL}{Deep Learning}
\acro{SGD}{Stochastic Gradient Descent}
\end{acronym}

\title{\LARGE \bf
Deep Forward-Backward SDEs for Min-max Control
}

\author{Ziyi Wang$^1$, Keuntaek Lee$^1$, Marcus A. Pereira$^1$, Ioannis Exarchos$^2$ and Evangelos A. Theodorou$^1$
\thanks{$^1$: Georgia Institute of Technology, Atlanta, GA, USA. $^2$ :Emory University, Atlanta, GA, USA. Email:zwang450@gatech.edu}
}

\begin{document}

\maketitle
\thispagestyle{empty}
\pagestyle{empty}

\begin{abstract}

This paper presents a novel  approach to numerically solve stochastic differential games  for nonlinear systems.  The proposed approach relies  on the nonlinear Feynman-Kac theorem that establishes a connection between parabolic deterministic partial differential equations and forward-backward stochastic differential equations. Using this theorem the Hamilton-Jacobi-Isaacs partial differential equation associated with differential games is represented by a system of forward-backward stochastic differential equations.  Numerical solution of the aforementioned system of stochastic differential equations is performed using importance sampling and a Long-Short Term Memory  recurrent neural network, which is trained in an offline fashion.  The resulting algorithm is tested on two example systems in simulation and compared against the standard risk neutral stochastic optimal control  formulations. 

\end{abstract}

\section{INTRODUCTION}
Stochastic optimal control is a mature discipline of control theory with a plethora  of applications to autonomy, robotics, aerospace systems, computational neuroscience, and finance.   From a methodological stand point, stochastic dynamic programming  is the pillar of stochastic optimal control theory. Application of the stochastic dynamic  programming  results in the so-called \ac{HJB} \ac{PDE}. Algorithms for stochastic control can be classified   into different categories depending on the way of how they are dealing with the curse of dimensionality in  solving the \ac{HJB} \ac{PDE} for systems with many degrees of freedom and/or states. 

Game-theoretic, or min-max, extension to optimal control was first investigated by Isaacs \cite{Isaacs1965}. He
associated the solution of a differential game with the solution to a HJB-like equation, namely its min-max extension, also known as the \ac{HJI} equation. The \ac{HJI} equation was derived heuristically under the assumptions of Lipschitz continuity of the cost and the dynamics, in addition to the assumption that both of them are separable in terms of the maximizing and minimizing controls. Despite extensive results in the theory of differential games, algorithmic development has seen less growth, due to the involved difficulties in addressing such problems. Prior work, including the Markov Chain approximation method \cite{Kushner2002}, largely suffers by the curse of dimensionality. In addition, a specific class of min-max control trajectory optimization methods have been derived recently, relying on the foundations of \textit{differential dynamic programming} (DDP)~\cite{Morimoto2003,Morimoto2002,Sun2015}, which requires linear and/or quadratic approximation of the dynamics and value function.

Due to the inherent difficulties of solving stochastic differential games, most of the effort in optimal control theory was focused on the \ac{HJB} \ac{PDE}. Addressing the solution of the \ac{HJB} equation, a number of algorithms for stochastic optimal control have been proposed that rely on the probabilistic representation of solutions of linear and nonlinear backward \acp{PDE}. Starting from the path integral control framework \cite{Kappen2005b}, the \ac{HJB} equation is transformed into a linear backward \ac{PDE} under certain conditions related to control authority and variance of noise. The  probabilistic representation of the solution of this \ac{PDE} is provided by the linear Feynman-Kac theorem \cite{Fleming1971,Fleming2006,Karatzasbook}. The nonlinear Feynman-Kac theorem avoids the assumption required in the path integral control framework at the cost, however, of representing the solution of the \ac{HJB} equation with a system of \acp{FBSDE} \cite{yong1999stochastic,Pardoux_Book2014}. Previous work by our group aimed at improving sampling efficiency and reducing computational complexity, and in \cite{exarchos2017stochastic,exarchos2016learning,exarchos2018stochastic} an importance sampling scheme was proposed and employed to develop  iterative stochastic control algorithms using the \ac{FBSDE} formulation. This work lead to algorithms for $L^{2}, L^{1} $, risk-sensitive stochastic optimal control, as well as stochastic differential games \cite{exarchosL1,exarchos2016game,exarchos2018DGA}. 

In \cite{han2016deep} the authors incorporate deep learning algorithms, such as Deep Feed-Forward Neural Networks, within the \ac{FBSDE} formulation and  demonstrated  the applicability the resulting algorithms to solving \acp{PDE}. While the approach in \cite{han2016deep}  offers an efficient method to represent the value function and its gradient, it has been only applied to \acp{PDE} that correspond to simple dynamics. Motivated by the limitations of the existing work on \acp{FBSDE} and \ac{DL}, the work in reference \cite{pereira2019neural} utilizes importance sampling together with the benefits of recurrent neural networks in order to capture the temporal dependencies of the value function and to scale the deep \ac{FBSDE} algorithm to high dimensional nonlinear stochastic systems. 

In this work, we demonstrate that the \acp{FBSDE} associated with stochastic differential games can be solved with the deep \ac{FBSDE} framework. We focus on the case of min-max stochastic control that corresponds to risk sensitive control. Using the \ac{LSTM} network architecture \cite{hochreiter1997lstm}, we introduce a scalable deep min-max \ac{FBSDE} controller that results in trajectories with reduced variance. We demonstrate the variance reduction benefit of this algorithm against the standard risk neutral stochastic optimal control formulation of the deep \ac{FBSDE} framework on a pendulum and a quadcopter in simulation.

The rest of this paper is organized as follows: in Section \ref{math} we introduce the min-max stochastic control problem, demonstrate its connection to risk sensitive control, and reformulate the problem with a system of \acp{FBSDE}. We present the min-max \ac{FBSDE} controller in Section \ref{algorithm}. In Section \ref{experiments}, we compare the controller introduced in this work against the deep \ac{FBSDE} algorithm for standard stochastic optimal control, and we explore the variance reduction benefit of our controller as a function of risk sensitivity. Finally, we conclude the paper in Section \ref{conclusions}.
   

\begin{figure}
  \includegraphics[width=\linewidth]{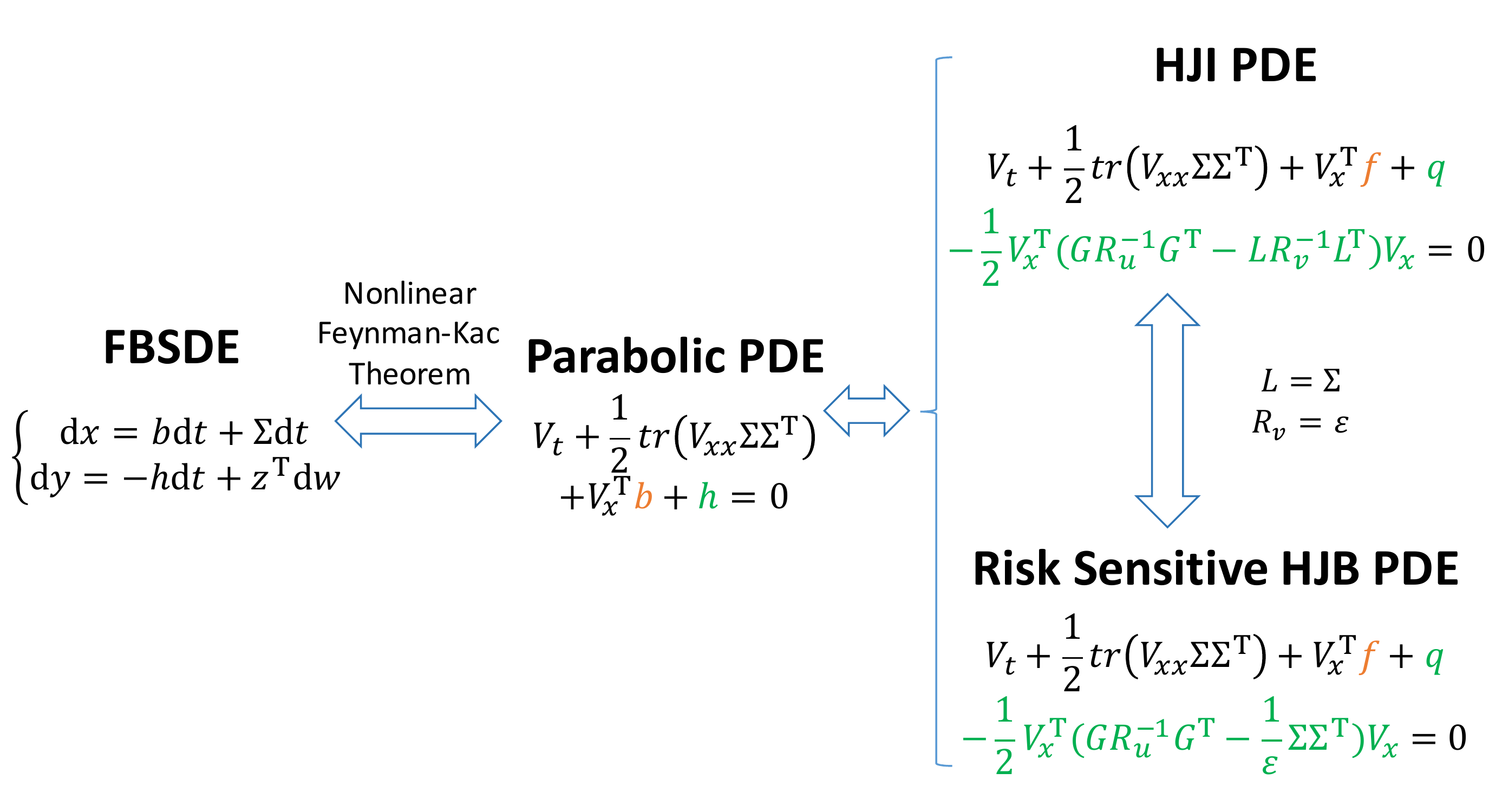}
  \caption{A schematic diagram showing the relationship between \acp{PDE} and \ac{FBSDE}. Terms in \textcolor{orange}{orange} denote drift in FSDE, and terms in \textcolor{green}{green} denote drift in BSDE.}
  \label{fig:schematic}
\end{figure}

\section{FBSDE for Differential Games}
\label{math}

\subsection{Min-Max Stochastic Control}
Consider a system with control affine dynamics in a differential game setting as follows:
\begin{align}\begin{split}
    \rd x &= f(x(t),t)\rd t + G(x(t),t) u(t) \rd t + L(x(t),t)v(t) \rd t \\
      & \hspace{3mm} + \Sigma(x(t),t) \rd w(t) \quad t  \in [\tau, T].
      \end{split}
\end{align}
where $\tau\in[0,T]$, $T$ is the task horizon, $x\in\R^n$ is the state, $u\in\R^p$ is the minimizing control, $v\in\R^q$ is the adversarial control, $w(t)$ is a standard $m$ dimensional Brownian motion, $f:\R^n\times[\tau, T]\rightarrow\R^n$ represents the drift, $G:\R^n\times[\tau, T]\rightarrow\R^{n\times p}$ represents the actuator dynamics, $L:\R^n\times [\tau, T]\rightarrow \R^{n\times q}$ represents the adversarial control dynamics and $\Sigma:\R^n\times[\tau, T]\rightarrow\R^{n\times m}$ represents the diffusion. For this system we can define the following cost function:
\begin{align}\begin{split} 
&J(\tau, x_{\tau}; u(\cdot), v(\cdot)) = \\
&\E \bigg[g(x(T))  + \int_{\tau}^{T} q(x(t),t)  + \cfrac{1}{2} u^{\rT} R_{u} u  - \cfrac{1}{2}v^{\rT} R_{v} v \rd t \bigg],
\end{split}
\end{align}
where $g:\R^n\rightarrow\R^+$ is the terminal state cost, $q:\R^n\rightarrow\R^+$ is the running state cost, and $R_u\in\R^{p\times p}$ and $R_v\in\R^{q\times q}$ are positive definite control cost matrices.

The min-max stochastic control problem is formulated as follows: 
\begin{equation}
    V(x_\tau,\tau) = \inf_{u \in \calU} \sup_{v \in \calV} J(\tau, x_{\tau}; u(\cdot), v(\cdot)),
\end{equation}
where the minimizing control's goal is reducing the cost under all admissible non-anticipating strategies $\calU$, while the adversarial control maximizes the cost under all admissible non-anticipating strategies $\calV$.

The \ac{HJI} equation for this problem is:  

\begin{equation}
      \begin{cases}
     V_t + \inf_{u \in \calU} \sup_{v \in \calV} \bigg\{\cfrac{1}{2}\tr\big(V_{xx} \Sigma \Sigma^\mathrm{T}\big) +  V_x^\mathrm{T}  ( f + G u + L v )     \\
    + q + \cfrac{1}{2} u^{\rT}  R_{u} u - \cfrac{1}{2} v^{\rT}  R_{v}  v \bigg\}  = 0,   \quad (t,x) \in  [\tau, T) \times \R^{n},  \\
    V(x,T) = g(x), \quad x \in \R^{n}.
    \end{cases}
\end{equation}

The terms inside the infimum and supremum operations are collectively called the Hamiltonian. The optimal minimizing and adversarial controls $ u $ and $ v $ are those for which the gradient of the Hamiltonian vanishes, which take the following form:
\begin{align}
    \begin{split}
        u(x(t),t) &= - R_{u}^{-1} G^{\rT} V_{x}, \\
         v(x(t),t) &=  R_{v}^{-1} L^{\rT} V_{x}. 
    \end{split}
    \label{eq:optimal_control}
\end{align}
 
Substitution of the expressions above into the \ac{HJI} equation results in:
\begin{equation}
      \begin{cases}
    V_t + \cfrac{1}{2}\tr\big(V_{xx}\Sigma \Sigma^\mathrm{T}\big)  
     - \cfrac{1}{2} V_{x}^{\rT} \bigg( G  R_{u}^{-1} G^{\rT} - L R_{v}^{-1} L^{\rT} \bigg) V_x \\ 
     +  V_x^\mathrm{T}  f + q = 0, \quad (t,x) \in  [\tau, T) \times \R^{n}, \\
    V(x,T)= g(x), \quad x \in \R^{n}.
    \end{cases}
    \label{eq:HJI_MinMax}
\end{equation}

Note that we will drop functional dependence in all \acp{PDE} for notational compactness. In the following section we show the equivalence of a certain case of min-max control to risk sensitive control.

\subsection{Risk Sensitive Stochastic Optimal Control}
Risk sensitive stochastic optimal control\cite{TamerBasar1995} is essential in cases where decision has to be made in a manner that is robust to the stochasticity of the environment.  Let us consider the following performance index: 
\begin{equation}
\begin{split}
    &J(\tau, x_{\tau}; u(\cdot)) =\\ &\epsilon \ln \E \bigg[\exp \cfrac{1}{\epsilon} \bigg( g(x(T)) + \int_{\tau}^{T} q(x(t),t) + \cfrac{1}{2} u(t)^{\rT} R  u(t) \rd t  \bigg) \bigg],
\end{split}
\end{equation}
where $\epsilon\in\R^+$ is the risk sensitivity. The risk sensitive stochastic optimal control problem is formulated with the following value function:
\begin{equation}
V(x_\tau, \tau)=\inf_{u \in \calU} J(\tau, x_{\tau}; u(\cdot)),
    \end{equation}
subject to the dynamics: 
\begin{equation}
\rd x(t) = f(x(t),t) \rd t + G(x(t),t) u(t) \rd t + \sqrt{\cfrac{\epsilon}{2 \gamma^{2}}} \tilde{\Sigma}(x(t),t) \rd w(t),
\end{equation}
where $\gamma\in\R^+$ is a small constant, and $\tilde{\Sigma}$ represents diffusion \cite{Fleming1995}.

The \ac{HJB} equation for this stochastic optimal control problem is formulated as follows:
\begin{equation}
      \begin{cases}
     V_t + \inf_{u\in\calU} \bigg\{ \cfrac{\epsilon}{4  \gamma^2}\tr\big(V_{xx} \tilde{\Sigma} \tilde{\Sigma}^\mathrm{T}\big) +  V_x^\mathrm{T}  ( f + G u )  + q   \\
    + \cfrac{1}{2} u^{\rT}  R u + \cfrac{1}{4 \gamma^{2}} V_x^{\rT}  \tilde{\Sigma} \tilde{\Sigma}^\mathrm{T} V_{x}  \bigg\}  = 0,   \quad (t,x) \in  [\tau, T) \times \R^{n},  \\
    V(x,T) = g(x), \quad x \in \R^{n}.
    \end{cases}
    \label{eq:HJB_rs_naive}
\end{equation}

The optimal control can be obtained by finding the control where the gradient of the terms inside the infimum vanishes and has the form $ u(x(t),t) = -R^{-1} G^{\rT} V_{x}$. By substituting in the optimal control and setting $ \Sigma = \sqrt{\cfrac{\epsilon}{2 \gamma^{2}}}\tilde{\Sigma}$ in \eqref{eq:HJB_rs_naive}, we get the following final form of the \ac{HJB} \ac{PDE}:  

\begin{equation}
      \begin{cases}
    V_t + \cfrac{1}{2}\tr\big(V_{xx}\Sigma \Sigma^\mathrm{T}\big)  
     - \cfrac{1}{2} V_{x}^{\rT} \bigg( G  R^{-1} G^{\rT} - \cfrac{1}{\epsilon} \Sigma \Sigma^{\rT} \bigg) V_x \\ 
     +  V_x^\mathrm{T}  f + q = 0, \quad (t,x) \in  [\tau, T) \times \R^{n}, \\
    V(x,T)= g(x), \quad x \in \R^{n}.
    \end{cases}\label{eq:HJB_Risk}
\end{equation}

Note that the above \ac{PDE} is a special case of the \ac{HJI} \ac{PDE} \eqref{eq:HJI_MinMax} when $ L = \Sigma $ and $ R_{v} = \epsilon I $ (Fig. \ref{fig:schematic}). Intuitively, this means that min-max control collapses to risk sensitive control when it is solving a problem with non-zero mean noise as the adversary, and the control authority of this adversary is proportional to the risk sensitivity.

\subsection{FBSDE Reformulation}

We now reformulate the min-max control \ac{PDE} \eqref{eq:HJI_MinMax} in the risk sensitive case to a set of \acp{FBSDE}. Here we restate the nonlinear Feynman-Kac theorem (Theorem 2) from \cite{exarchos2016game}:

\begin{theorem}[\text{Nonlinear Feynman-Kac}]
\label{nfk}
\textit{Consider the following Cauchy problem:}
\begin{equation}
\begin{cases}
    V_t + \cfrac{1}{2}\tr(V_{xx}\Sigma\Sigma^\rT) + V_x^\rT b + h=0, \ \ 
    (t,x) \in [\tau, T)\times \R^n,\\ V(T,x) = g(x), \quad x\in \R^n,
\end{cases}
\label{eq:parabolic_pde}
\end{equation}
\textit{wherein the functions $\Sigma,\, b(t,x), \, h(t,x,V,\Sigma^\rT V_x)$, and $g(x)$ satisfy mild regularity conditions. Then \eqref{eq:parabolic_pde} admits a unique viscosity solution $V:[\tau, T]\times \R^n \rightarrow \R$, which has the following probabilistic representation:}
\begin{equation}
\Big(V(t,x),\, \Sigma^\rT V_x(t,x) \Big) = \Big(y(t,x),\, z(t,x)\Big), \ \forall (t,x)\in [\tau, T]\times \R^n,
\end{equation}
\textit{wherein $\big(x(\cdot), y(\cdot), z(\cdot)\big)$ is the unique adapted solution of the \acp{FBSDE} given by:}
\begin{equation}
\begin{cases}
    \rd x(t) = b(x(t),t)\rd t + \Sigma(x(t),t)\rd w(t), \quad t\in[\tau, T] \\
    x(\tau) = \xi
\end{cases}
\label{eq:fsde}
\end{equation}
\textit{and}
\begin{equation}
\begin{cases}
    \rd y(t) = -h(t, x(t), y(t), z(t))\rd t + z(t)^\rT \rd w(t), \quad t\in[\tau, T] \\
    y(T) = g(x(T))
\end{cases}
\label{eq:bsde}
\end{equation}
\end{theorem}
In order to apply the Nonlinear Feynman-Kac theorem to \eqref{eq:HJI_MinMax}, we assume that there exist matrix-valued functions $\Gamma_u: [\tau, T]\times\R^n\rightarrow \R^{m\times p}$ and $\Gamma_v: [\tau, T]\times\R^n\rightarrow \R^{m\times q}$ such that $G(x(t),t)=\Sigma(x(t),t)\Gamma_u(x(t),t)$ and $L(x(t),t)=\Sigma(x(t),t)\Gamma_v(x(t),t)$ for all $(t,x)\in[\tau, T]\times \R^n$, satisfying the same regularity conditions. This assumption suggests that there can not be a channel containing control input but no noise. In the risk sensitive case of min-max control, this assumption is already satisfied with $L(x(t),t)=\Sigma(x(t),t)$ and $\Gamma_v(x(t),t)=I$, where $I$ is a $m \times m$ identity matrix because adversarial control enters through the noise channels. Under this assumption, Theorem \ref{nfk} can be applied to the risk sensitive case of \ac{HJI} equation \eqref{eq:HJI_MinMax} with
\begin{align}
    \begin{split}
        b(x(t),t) &= f(x(t),t)\\
        h(x(t), y(t), z(t), t) &= q(x(t))\\ &- \cfrac{1}{2}V_x^\rT\Big(\Sigma\Gamma_u R^{-1}\Gamma_u^\rT\Sigma^\rT - \cfrac{1}{\epsilon}\Sigma\Sigma^\rT\Big)V_x.
    \end{split}
    \label{eq:rs_FBSDE_coeffs}
\end{align}
The relationship between \ac{FBSDE} \eqref{eq:fsde}, \eqref{eq:bsde}, HJI PDE \eqref{eq:HJI_MinMax}, HJB PDE \eqref{eq:HJB_Risk}, and the parabolic PDE \eqref{eq:parabolic_pde} is summarized in Fig. \ref{fig:schematic}.

\subsection{Importance Sampling}
The system of \acp{FBSDE} in \eqref{eq:fsde} and \eqref{eq:bsde} corresponds to a system whose dynamics are uncontrolled. In many cases, especially for unstable systems, it is hard or impossible to reach the target state with uncontrolled dynamics. We can address this problem by modifying the drift term in the dynamics (forward SDE) with an additional control term. Through Girsanov's theorem \cite{Girsanov} on change of measure, the drift term in the forward SDE \eqref{eq:fsde} can be changed if the backward SDE \eqref{eq:bsde} is compensated accordingly. This results in a new \ac{FBSDE} system given by
\begin{align}
    \begin{cases}
    \rd \tilde{x}(t) &= [b(\tilde{x}(t),t) + \Sigma(\tilde{x},t)K(t)]\rd t \\ & + \Sigma(\tilde{x}(t),t)\rd \tilde{w}(t), \quad t\in[\tau, T] \\
    \tilde{x}(\tau) &= \xi
    \end{cases}
    \label{eq:new_fsde}
\end{align}
and
\begin{align}
    \begin{cases}
    \rd \tilde{y}(t) &= [-h(t, \tilde{x}(t), \tilde{y}(t), \tilde{z}(t)) + \tilde{z}^\rT K(t)]\rd t \\ & + \tilde{z}(t)^\rT \rd \tilde{w}(t), \quad t\in[\tau, T] \\
    \tilde{y}(T) &= g(\tilde{x}(T))
    \end{cases},
    \label{eq:new_bsde}
\end{align}
for any measurable, bounded and adapted process $K:[\tau, T] \rightarrow \R^n$. It is easy to verify that the \ac{PDE} associated with the new system is the same as the original one \eqref{eq:parabolic_pde}. For the full derivation of change of measure for \acp{FBSDE}, we refer readers to proof of Theorem 1 in \cite{exarchos2018stochastic}. We can conveniently set $K=\Gamma_u(\tilde{x}(t), t)\bar{u} + \Gamma_v(\tilde{x}(t), t)\bar{v}$ for min-max control. Note that the nominal controls $\bar{u}$ and $\bar{v}$ can be any open or closed loop control or control from a previous iteration.

\begin{figure*}[h]
\centering
  \includegraphics[width=450pt]{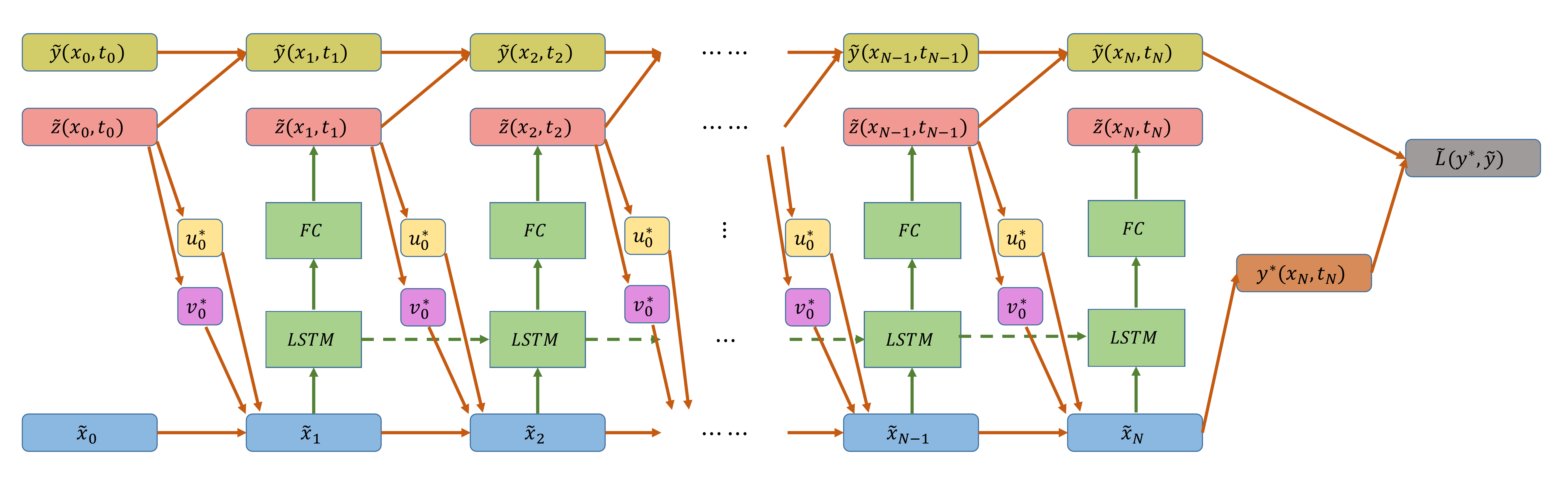}
  \caption{LSTM network architecture.}
\label{fig:LSTMnetwork}
\end{figure*}

\begin{table*}[!ht]
\centering
\caption{Comparison of total state variance between deep min-max controller and baseline deep FBSDE controller.}
\label{table:1}
\begin{tabular}{ |l||c|c|c|c|  }
 \hline
  & \multicolumn{2}{c|}{Pendulum} & \multicolumn{2}{c|}{QuadCopter} \\
 \hline
  & Low Noise & High Noise & Low Noise & High Noise \\
 \hline
 Baseline   & 5.3 & 149.5 & 3.2 & 78.6\\
 RS         & 3.9 & 134.2 & 2.7 & 69.6\\
 \hline
 Variance Reduction (\%) & 26 & 10 & 16 & 11\\
 \hline
\end{tabular}
\end{table*}

\subsection{Forward Sampling of BSDE}
The compensated BSDE \eqref{eq:new_bsde} needs to satisfy a terminal condition, meaning its solution needs to be propagated backward in time, yet the filtration evolves forward in time. This poses a challenge for sampling based methods to solve the system of \acp{FBSDE}. One solution is to approximate the conditional probability of the process and back-propagate the expected value. This approach lacks scalability due to inevitable compounding of approximation errors that are accumulated at every time step during regression.

This problem can be alleviated with \ac{DL}. Using a deep recurrent network, we can initialize the value function and its gradient at $t=\tau$ and treat the initializations as trainable network parameters. This allows for the BSDE to be propagated forward in time along with the FSDE. At the final time, the terminal condition can be compared against the propagated value in the loss function to update the initialzations as well as the network parameters. Compared to the conditional probability approximation scheme, the \ac{DL} approach has the additional advantage of not accumulating errors at every time step since the recurrent network at each time step contributes to a common goal of predicting the target terminal condition and thus prediction errors are jointly minimized.

\section{Deep Min-max FBSDE Controller}
\label{algorithm}

With eqs. \eqref{eq:new_fsde} and \eqref{eq:new_bsde}, we have a system of \acp{FBSDE} that we can sample from around a nominal control trajectory. Inspired by the network architecture developed in \cite{pereira2019neural}, we propose a deep min-max FBSDE algorithm that solves the risk sensitive formulation of the min-max control.

\subsection{Numerics and Network Architecture}
The task horizon $\tau <t<T$ can be discretized as $n=\{0, 1, \cdots, N\}$ with a time discretization of $\Delta t = (T-\tau)/N$. With this we can approximate the continuous variables as step functions and obtain their discretization as $\tilde{x}_n, \tilde{y}_n, \tilde{z}_n, u_n = \tilde{x}(t), \tilde{y}(t), \tilde{z}(t), u(t)$ if $ \tau + n\Delta t \leq t < \tau + (n+1)\Delta t$.

The network architecture used in this paper is shown in Fig. \ref{fig:LSTMnetwork}, which is based on the \ac{LSTM} network in \cite{pereira2019neural} with min-max dynamics and value function dynamics incorporated. \ac{LSTM} is a natural choice of network here since it is designed to effectively deal with the vanishing gradient problem in recurrent prediction of long time series \cite{hochreiter1997lstm}. We use a two-layer \ac{LSTM} network with tanh activation and Xavier initialization \cite{glorot2010understanding}. At every time step, the \ac{LSTM} predicts the value function gradient using the current state as input. The optimal minimizing and adversarial control are then calculated with
\begin{align}
    u_n^* &= -R_u^{-1}\Gamma_{u,n}\tilde{z}_n\\
    v_n^* &= \cfrac{1}{\epsilon}\tilde{z}_n
\end{align}
and fed back to the dynamics for importance sampling. Note that the adversarial control is only present during training. After the network is trained, only the optimal minimizing control is used at test time. By exposing the minimizing controller to an adversary that behaves in an optimal fashion, it becomes more robust resulting in trajectories with smaller variances.

\subsection{Algorithm}

\begin{algorithm}[h]
\caption{Deep Min-max FBSDE Controller}
\begin{algorithmic}
 \STATE \textbf{Given:}\\ 
          $\tilde{x}_0=\xi,\;f,\; G, \Sigma $: Initial state and system dynamics;\\ $g,\; q, \; R_u, \; \epsilon(=R_v)$: Cost function parameters;\\$N$: Task horizon; $K$: Number of iterations; $M$: Batch size; \\ $\Delta t$: Time discretization;
          $\lambda$: weight-decay parameter;\\ $\gamma$: Loss function parameter;
          \STATE \textbf{Parameters:}\\
          $\tilde{y}_0=V(\tilde{x}_0,\tau;\psi)$: Value function at $t=\tau$;\\
          $\tilde{z}_0=\Sigma^\mathrm{T} \,\nabla_{\tilde{x}} V $: Gradient of value function at $t=\tau$;\\
          $\theta$: Weights and biases of all LSTM layers;
          \STATE \textbf{Initialize:}\\ $\{\tilde{x}^i_0\}_{i=1}^{M},\;\tilde{x}^i_0=\xi $\\
          $\{\tilde{y}^i_0\}_{i=1}^{M},\;\tilde{y}^i_0=V(\tilde{x}^i_0,0;\psi)$\\
          $\{\tilde{z}^i _0\}_{i=1}^{M},\;\tilde{z}^i _0=\Sigma^\mathrm{T} \nabla_{\tilde{x}} V(\tilde{x}^i_0,0;\psi)$
          \FOR{$k=1$ \TO $K$}
            \FOR{$i=1$ \TO $M$}
                \FOR{$n=1$ \TO $N-1$}
                    \STATE Compute gamma matrix: $\Gamma_{u,n}^i=\Gamma_u \big(\tilde{x}^i_n\big)$; 
                    \STATE $u^{i*}_n=-R_u^{-1}\Gamma^{i\mathrm{T}}_{u,n} \tilde{z}^i_n$;
                    \STATE $v^{i*}_n=\cfrac{1}{\epsilon} \tilde{z}^i_n$;
                    \STATE Sample Brownian noise: $\Delta \tilde{w}^i_n \sim \mathcal{N}(0, 1)$
                    \STATE Update value function:
                    \STATE $\tilde{y}^i_{n+1} = \tilde{y}^i_n - \Big( \tilde{h}\big(\tilde{x}^i _n, \, \tilde{y}^i_n,\,\tilde{z}^i_n\big) + \tilde{z}^i_n (\Gamma^i_{u,n} u^{i*}_n + v^{i*}_n) \Big) \Delta t $
                    \STATE $+ \tilde{z}^{i\mathrm{T}}_n \Delta \tilde{w}^i_n \sqrt{\Delta t} $
                    \STATE Update system state:
                    \STATE $\tilde{x}^i_{n+1}=\tilde{x}^i_n + f(\tilde{x}^i_n) \Delta t + \Sigma \big((\Gamma^i_{u,n} u^{i*}_n +  v^{i*}_n)\Delta t +\Delta \tilde{w}^i_n\sqrt{\Delta t}\big) $
                    \STATE Predict gradient of value function: $\tilde{z}^i_{n+1} = f_{LSTM}\big(\tilde{x}^i_{n+1}; \theta_k\big)$
                \ENDFOR
                \STATE Compute target terminal value: $y^{*i}_N=g\big(\tilde{x}^i_N\big)$
            \ENDFOR
            \STATE Compute mini-batch loss: $\tilde{\mathcal{L}}= \displaystyle \cfrac{1}{M} \sum_{i=1}^M \Big( \gamma\| y^{*i}_N - \tilde{y}^i_N\|^2 + (1-\gamma)\|y^{*i}_N\|^2 \Big) + \lambda \|\theta_k^2 \|^2 $
            \STATE $\theta_{k+1}\leftarrow$ Adam.step($\mathcal{L}, \theta_k$); $\psi_{k+1}\leftarrow$ Adam.step($\mathcal{L},\psi_k$)
          \ENDFOR
          \RETURN $\theta_K, \psi_K$
\end{algorithmic}
\label{alg:FBSDEalgorithm}
\end{algorithm}

The Deep Min-max FBSDE algorithm can be found in Algorithm \ref{alg:FBSDEalgorithm}. It solves a finite time horizon control problem by approximating the gradient of the value function $\tilde{z}^i_n$ (the superscript $i$ denotes the batch index, and the batch-wise computation can be done in parallel) at every time step with a \ac{LSTM}, which is parameterized by $\theta$, and propagating the \ac{FBSDE} associated with the control problem. For a given initial state condition $\xi$, the algorithm randomly initializes the value function and its gradient at $n=0$. The initial values are trainable and are parameterized by $\psi$. During training, at every time step, control inputs are sampled around the optimal minimizing and adversarial controls and applied to the system. The discretized forward dynamics and the value function SDEs are propagated using an explicit forward Euler integration scheme. The function $h$ is calculated using \eqref{eq:rs_FBSDE_coeffs}. At the final time step $n=N$, a modified $L^2$ loss with regularization is computed which compares the propagated value function $\tilde{y}_N^i$ against the true value function $y^{*i}_N$ calculated using the final state ($y^{*i}_N=g(\tilde{x}_N^i)$). For training our network, we propose a new regularized loss function, which is a convex combination of a) the difference between the target and the predicted value function, and b) the target value function itself:
\begin{equation}
    \tilde{\mathcal{L}}= \displaystyle \cfrac{1}{M} \sum_{i=1}^M \Big( \beta \| y^{*i}_N - \tilde{y}^i_N\|^2 + (1-\beta)\|y^{*i}_N\|^2 \Big) + \text{Reg}(\theta_k),
\end{equation}
since we want the prediction to be close to the target and at the same time, the target value function to converge to zero for the sake of the optimality. Notice that this additional component in the loss function is possible only due to importance sampling. The modified drift is implemented as a connection in the computational graph between the LSTM output and input to forward SDE at the next timestep. This allows the network parameters to influence the next state and hence the final state. The network can be trained by \ac{SGD} type optimizer and in our experiments, we used the Adam \cite{Adam} optimizer.


\section{Experiments}
\label{experiments}

The algorithm is implemented on a pendulum and quadcopter system in simulation. The task for the two systems is to reach a target state. The trained networks are tested on 128 trajectories. The time discretization is 0.02 seconds across all cases. We compare the algorithm proposed in this paper with the one in \cite{pereira2019neural}, where the standard optimal control problem is considered, in two different noise conditions. We will use ``RS'' to denote the algorithm in this work and ``Baseline'' for the algorithm that we are comparing against. All experiments were done in TensorFlow \cite{tensorflow} on an Intel i7-4820k CPU Processor.

\begin{figure}
  \includegraphics[width=\linewidth]{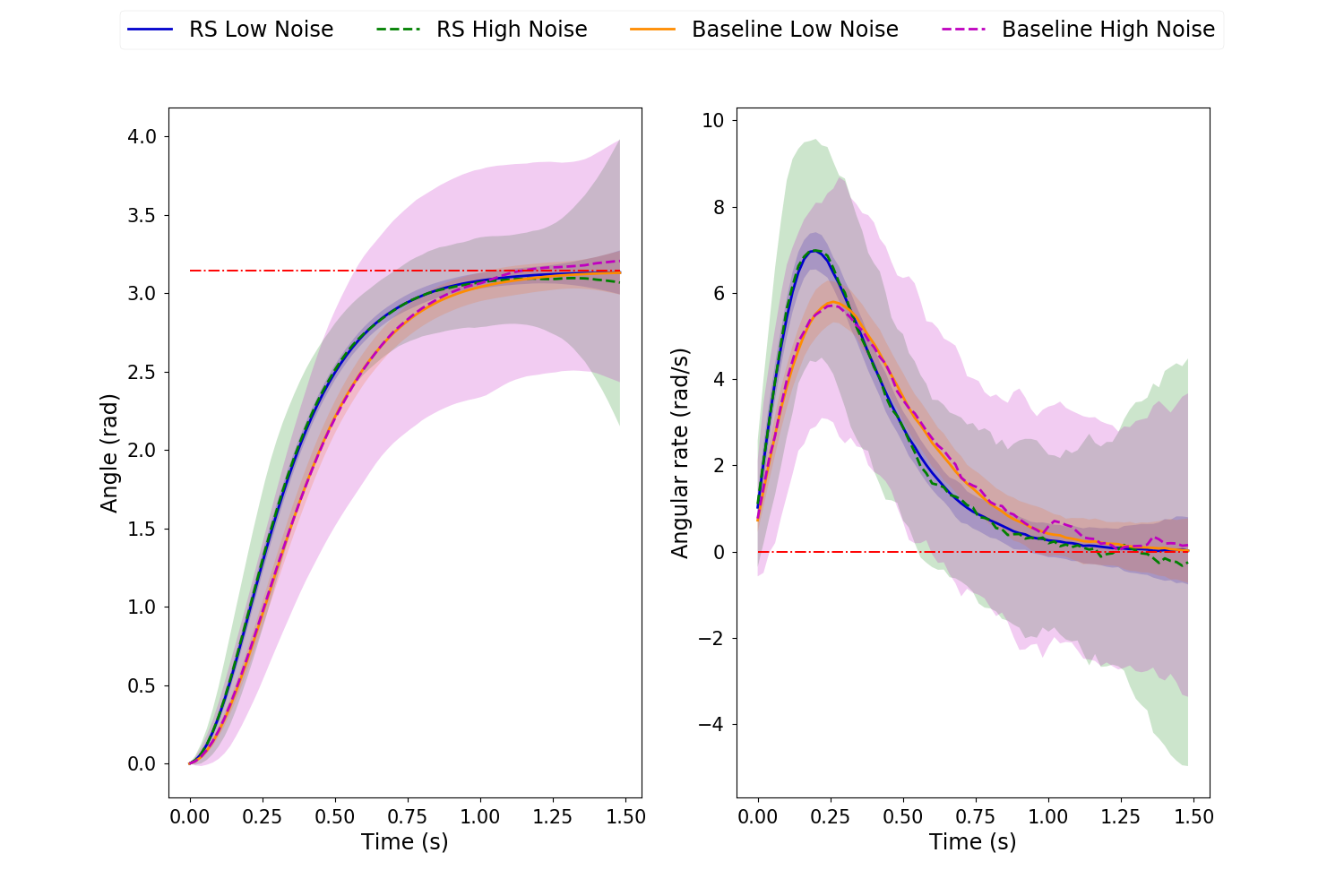}
  \caption{Pendulum states. \textit{Left}: Pendulum Angle; \textit{Right}: Pendulum Angular Rate.}
  \label{fig:pendulumStates}
\end{figure}
\begin{figure}
  \includegraphics[width=\linewidth]{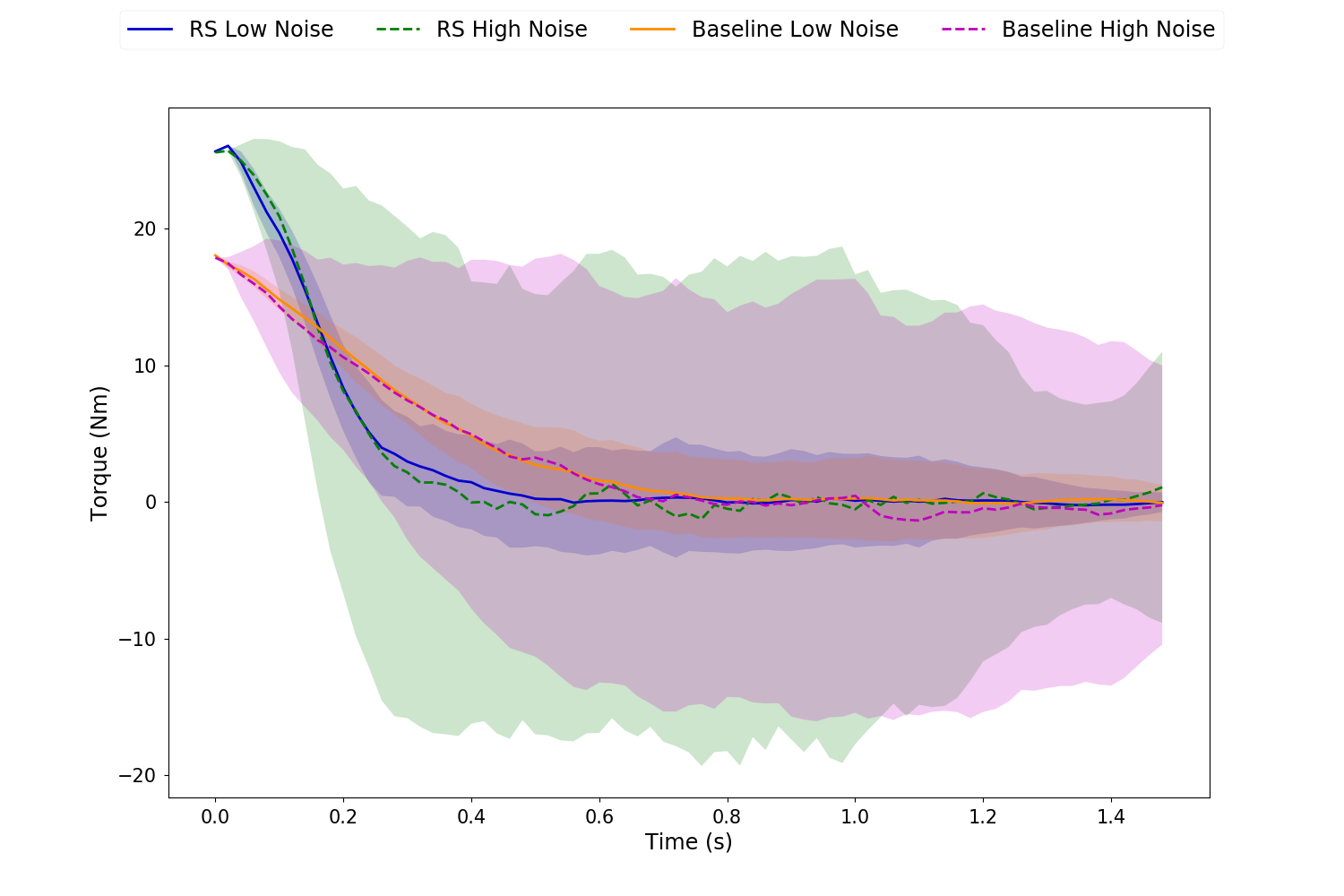}
  \caption{Pendulum controls.}
  \label{fig:pendulumControls}
\end{figure}

In all trajectory plots, the solid line denotes the mean trajectory in low noise condition, the dashed line denotes the mean trajectory in high noise condition, and the \textcolor{red}{red} dashed line denotes the target state. In addition, the 4 conditions are denoted by different colors, with \textcolor{blue}{blue} for RS in low noise condition, \textcolor{green}{green} for RS in high noise condition, \textcolor{orange}{orange} for Baseline in low noise condition, and \textcolor{magenta}{magenta} for Baseline with high noise. The shaded region of each color denotes the 95\% confidence region.

\subsection{Pendulum}
For the pendulum system, the algorithm was implemented to complete a swing-up task with a task horizon of 1.5 seconds. The two system states are the pendulum angle [$rad$] and the pendulum angular rate [$rad/s$]. Fig. \ref{fig:pendulumStates} plots the pendulum states in all 4 cases (RS with low and high noise and Baseline with low and high noise). The control applied to the system is the torque [$N \cdot m$] (Fig. \ref{fig:pendulumControls}).

\begin{figure}
  \includegraphics[width=\linewidth]{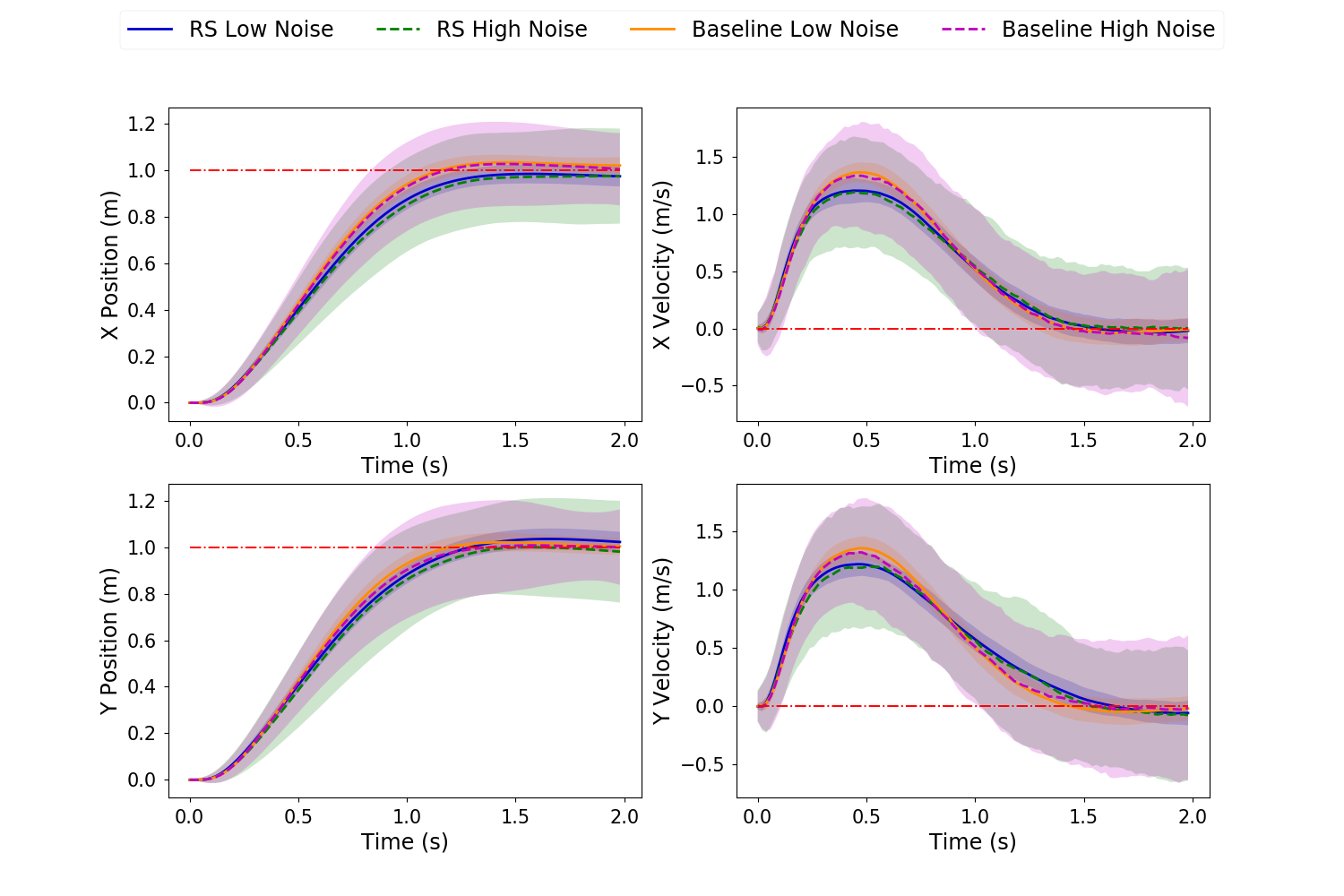}
  \caption{Quadcopter states. \textit{Top Left}: X Position; \textit{Top Right}: X Velocity; \textit{Bottom Left}: Y Position; \textit{Bottom Right}: Y Velocity.}
  \label{fig:QuadStates1}
\end{figure}
\begin{figure}
  \includegraphics[width=\linewidth]{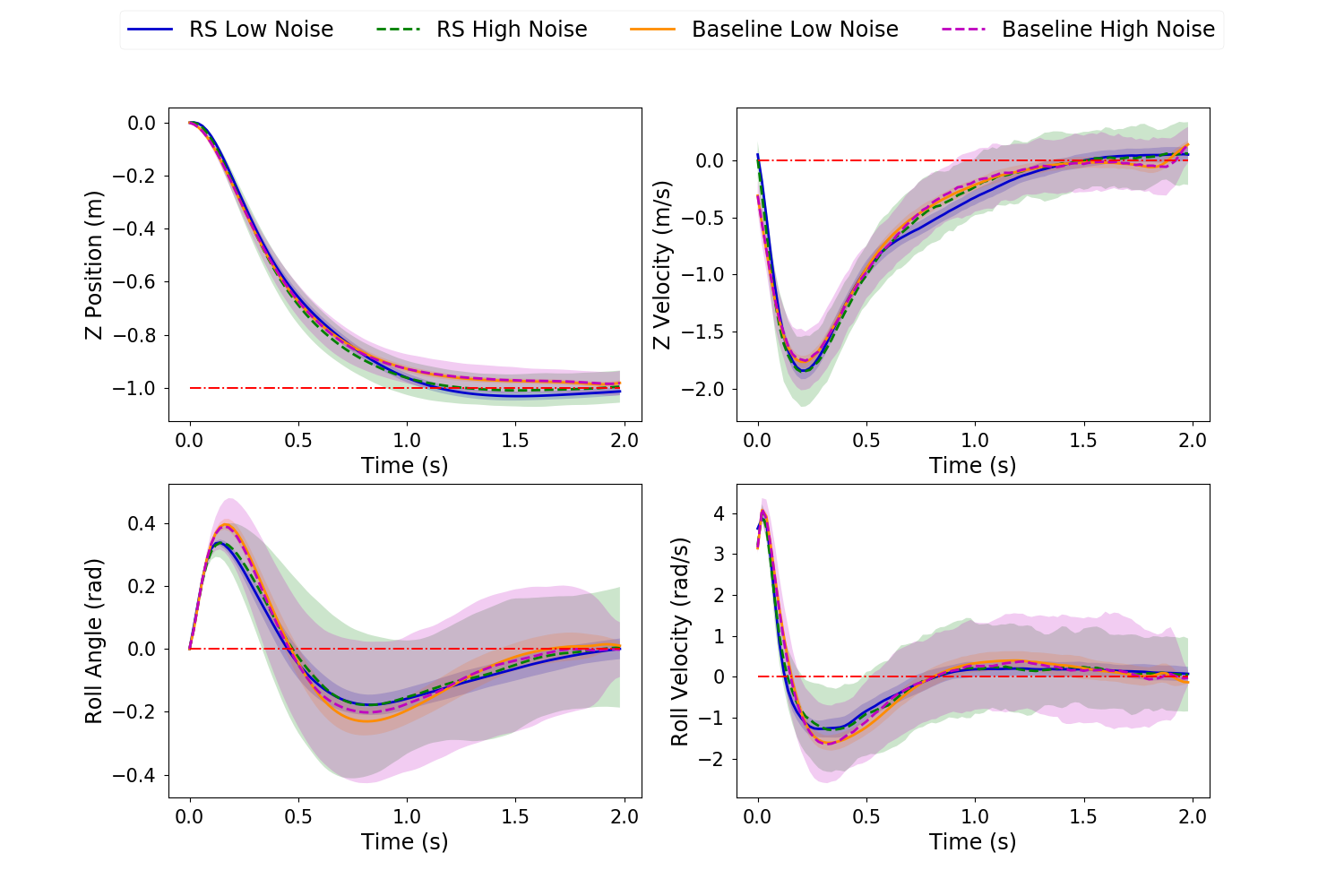}
  \caption{Quadcopter states. \textit{Top Left}: Z Position; \textit{Top Right}: Z Velocity; \textit{Bottom Left}: Roll Angle; \textit{Bottom Right}: Roll Velocity.}
  \label{fig:QuadStates2}
\end{figure}
\begin{figure}
  \includegraphics[width=\linewidth]{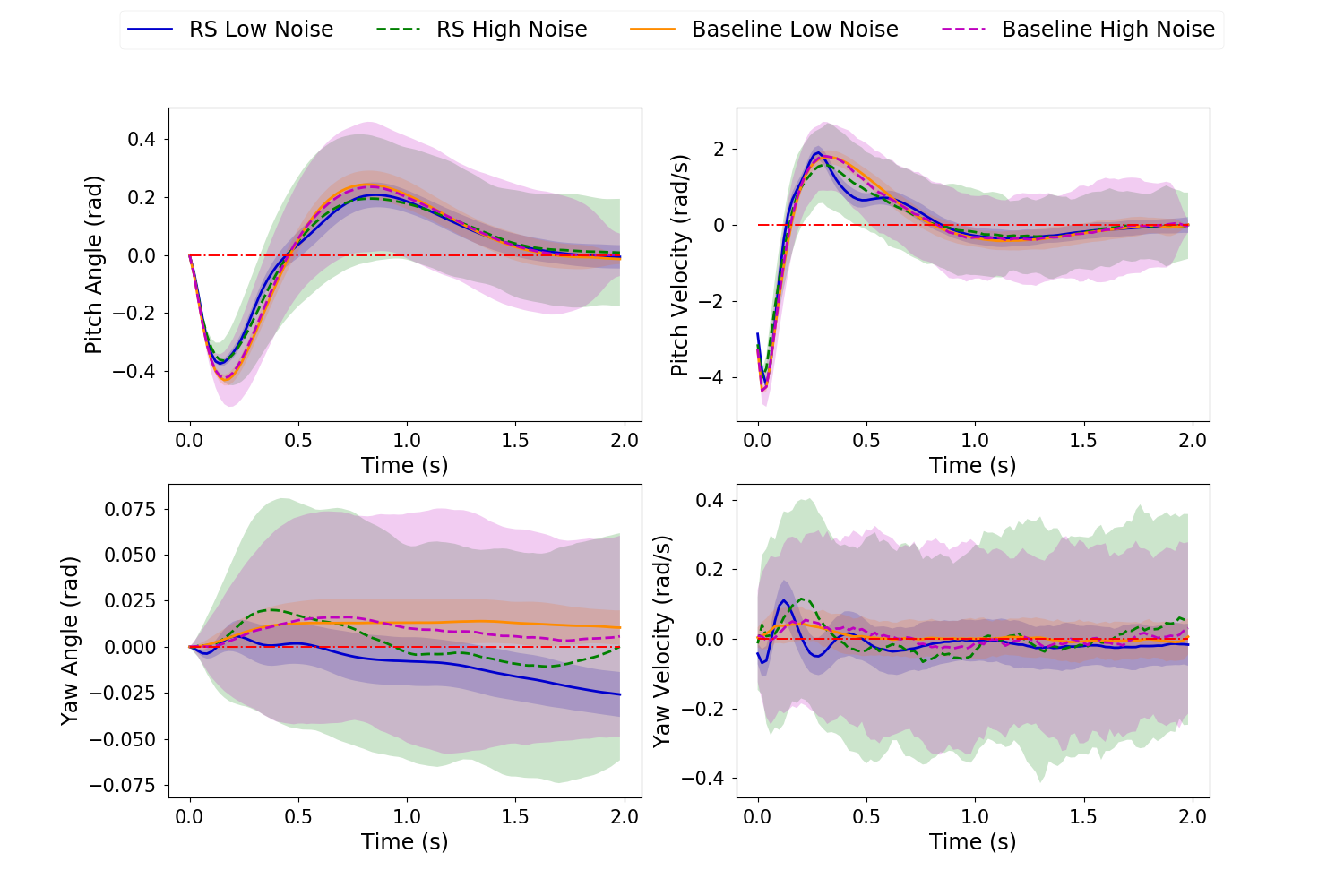}
  \caption{Quadcopter states. \textit{Top Left}: Pitch Angle; \textit{Top Right}: Pitch Velocity; \textit{Bottom Left}: Yaw Angle; \textit{Bottom Right}: Yaw Velocity.}
  \label{fig:QuadStates3}
\end{figure}
\begin{figure}
  \includegraphics[width=\linewidth]{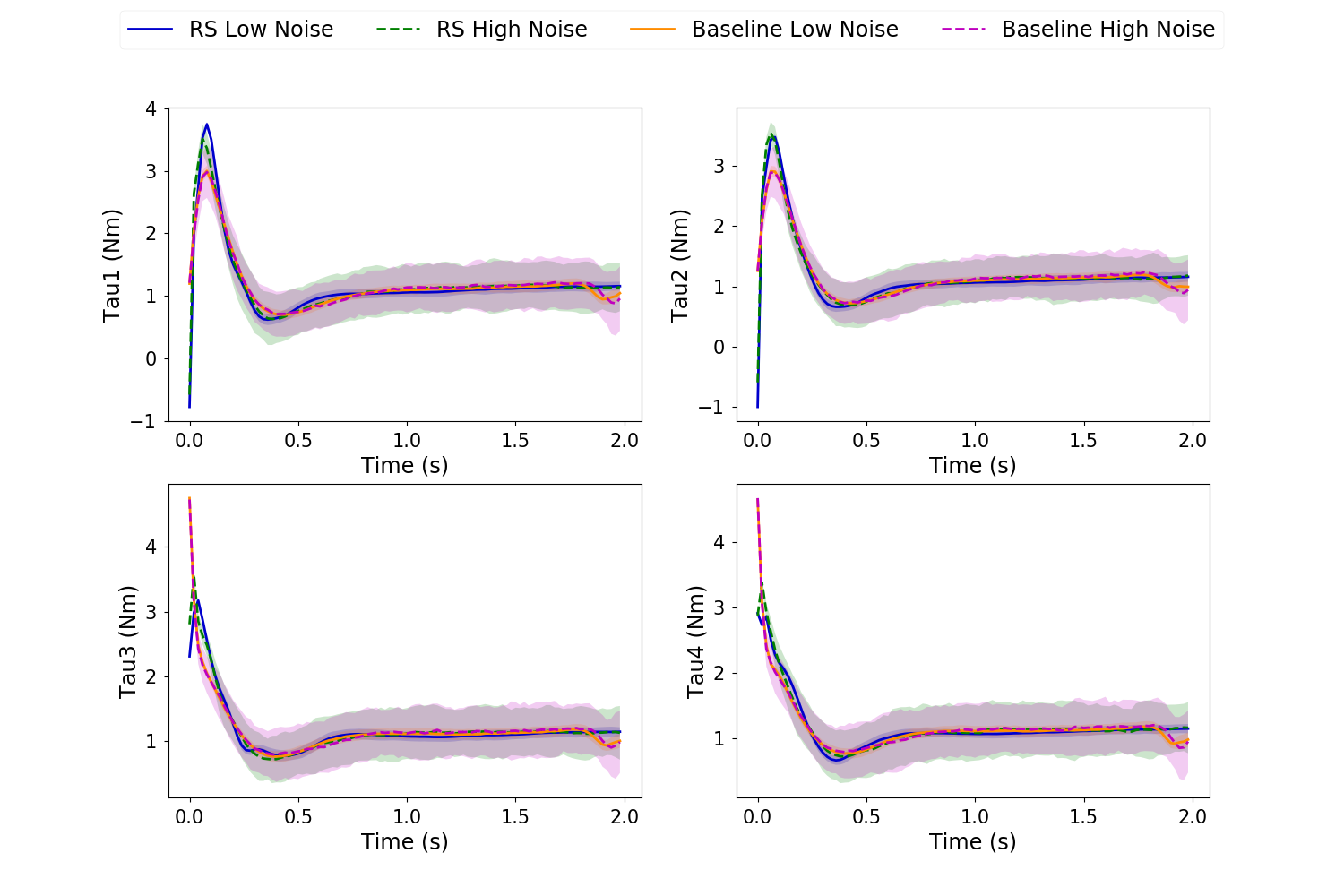}
  \caption{Quadcopter controls.}
  \label{fig:QuadControls}
\end{figure}
\begin{figure*}[ht]
\centering
  \includegraphics[width=\linewidth]{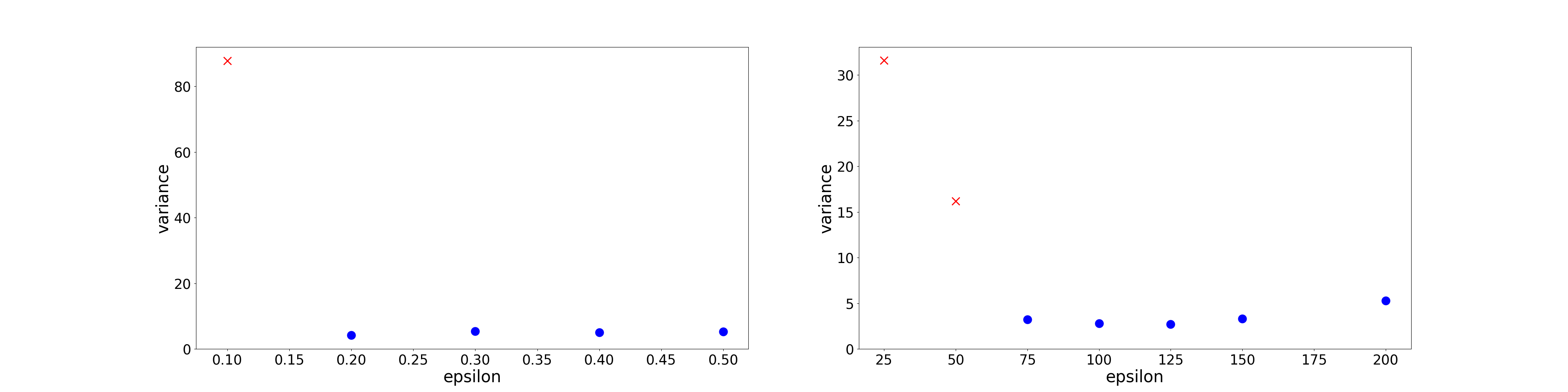}
  \caption{Total state variance vs. $\epsilon$ for both systems. \textit{Left: Pendulum; Right: QuadCopter.}}
  \label{fig:Variance}
\end{figure*}


\subsection{Quadcopter}
The algorithm was implemented on a quadcopter system for the task of reaching a final target state from an initial position with a task horizon of 2 seconds. The initial condition is 0 across all states. The target is 1 [$m$] upward, forward and to the right from the initial position with zero velocities and attitudes. The quadcopter dynamics used can be found in \cite{quad_dynamics}. The 12 system states are composed of the position [$m$], angles [$rad$], linear velocities [$m/s$], and angular velocities [$rad/s$]. The control inputs to the system are 4 torques [$N \cdot m$], which control the rotors (Fig. \ref{fig:QuadControls}).

\subsection{Reduced Variance with Deep Min-max FBSDE Controller}
The trajectory plots (Fig. \ref{fig:pendulumStates}, \ref{fig:QuadStates1}, \ref{fig:QuadStates2}, and \ref{fig:QuadStates3}) compare the Deep Min-max \ac{FBSDE} controller against the risk neutral Deep \ac{FBSDE} controller in a low noise setting and a high noise setting for both systems. From the plots we can observe that the min-max controller proposed in this work accomplishes the tasks with similar level of performance compared to the baseline controller. Numerical comparisons of the total state variance (sum of variance in all states over the entire trajectory) of all test cases can be found in Table \ref{table:1}. The results demonstrate at least 10\% reduction in total state variance across all cases. It is worth noting that the high noise setting results in less variance reduction benefits. By examining the substitution of $\Sigma=\sqrt{\cfrac{\epsilon}{2\gamma^2}}\tilde{\Sigma}$ from \eqref{eq:HJB_rs_naive} to \eqref{eq:HJB_Risk} in risk sensitive control derivation, we can see that increasing noise level is in some sense equivalent to increasing $\epsilon$. This naturally reduces the effect of the risk sensitive controller, as shown in the next section.

\subsection{Variance vs. Risk Sensitivity}
We also investigated the relationship between total state variance and risk sensitivity in the two systems. Fig. \ref{fig:Variance} plots the total state variance for different $\epsilon$ $(R_v)$ values while also keeping track of task completion. In the variance versus $\epsilon$ scatter plots, \textcolor{blue}{blue circles} are used to denote runs with successful task completion, whereas \textcolor{red}{red cross} denotes runs where the task failed. Since the risk sensitivity parameter $\epsilon$ is inversely proportional to the adversarial control authority, we expect the risk sensitive controller to converge to standard optimal controller as $\epsilon$ increases to infinity. On the other hand, as $\epsilon$ gets smaller, the adversarial control will eventually dominate the minimizing control and cause task failure. This is reflected in the plots as we can observe that the minimizing controller starts to fail when $\epsilon$ is too low. It is worth noting that the failure threshold increases as the system gets more complex and higher dimensional. Although the variance starts to increase as $\epsilon$ increases in the Quadcopter plot, the convergence to standard optimal controller is harder to observe as we only explore a limited range of $\epsilon$ values.



\section{Conclusions}
\label{conclusions}
In this paper, we proposed the Deep Min-max \ac{FBSDE} Control algorithm, based on the risk sensitive case of stochastic game-theoretic optimal control theory. Utilizing prior work on importance sampling of \acp{FBSDE} and efficiency of the \ac{LSTM} network to predict long time series, the algorithm is capable of solving stochastic game-theoretic control problems for nonlinear systems with control-affine dynamics. Comparison of this algorithm against the standard stochastic optimal control formulation suggests that by considering an adversarial control in the form of noise-related risk, the controller outputs trajectories with lower variance. Our algorithm scales in terms of the number of system states and system complexity for the min-max control problem, while the previous works did not. For future works, we would like to explore different network architectures to reduce the training time.



\addtolength{\textheight}{-3cm}   

\bibliographystyle{unsrt}
\bibliography{references}

\end{document}